\documentclass[journal]{IEEEtran}
\ifCLASSINFOpdf
\else
   \usepackage[dvips]{graphicx}
\fi
\usepackage{url}
\usepackage{amsmath}
\newtheorem{thm}{Theorem}
\newtheorem{cor}{Corollary}
\newtheorem{lem}{Lemma}

\usepackage{cite}
\usepackage{graphicx}
\usepackage{diagbox}
\usepackage{tabularx}
\usepackage{booktabs}
\usepackage{color}
\usepackage{epstopdf}

\begin{document}

\title{New-Type Hoeffding's Inequalities and Application in Tail Bounds }

\author{Pingyi Fan, \IEEEmembership{Senior Member, IEEE},
\thanks{Pingyi Fan are with Beijing National Research Center for Information Science and Technology and the Department of Electronic Engineering, Tsinghua University, Beijing 10084, China(e-mail:fpy@tsinghua.edu.cn).}
}

\maketitle

\begin{abstract}
It is well known that Hoeffding's inequality has a lot of applications in the signal and information processing fields. How to improve Hoeffding's inequality and find the refinements of its applications have always attracted much attentions. An improvement of Hoeffding inequality was recently given by Hertz \cite{r1}. Eventhough such an improvement is not so big, it still can be used to update  many known results with original Hoeffding's inequality, especially for Hoeffding-Azuma inequality for martingales. However, the results in original Hoeffding's inequality and its refinement one by Hertz only considered the first order moment of random variables. In this paper, we present a new type of Hoeffding's inequalities, where the high order moments of random variables are taken into account.  It can get some considerable improvements in the tail bounds evaluation compared with the known results. It is expected that the developed new type Hoeffding's inequalities could get more interesting applications in some related fields that use Hoeffding's results.
\end{abstract}

\begin{IEEEkeywords}
Hoeffding's Lemma, Hoeffding's tail bounds, Azuma inequality, Chernoff's bound.
\end{IEEEkeywords}

\IEEEpeerreviewmaketitle

\section{Introduction}

It is well known that Hoeffding's inequality has been applied in many scenarios in the signal and information processing fields.
Since Hoeffding's inequality was first found in 1963 \cite{r2}, it has been attracting much attentions in the academic research \cite{r3} \cite{r4} and industry.  Especially, in the last decade, it has been used to evaluate the channel code design \cite{r5}\cite{r6} and achievable rate over nonlinear channels \cite{r7} as well as delay performance in CSMA with linear virtual channels under a general topology \cite{r8} in information theory \cite{r9}. As one key tool, it also found the applications in machine learning and big data processing, i.e. PAC-Bayesian method analysis and Markov model analysis in machine learning \cite{r10}\cite{r11},statistical mode bias analysis\cite{r12}, concept drift in online learning for big data mining \cite{r13} and compressed sensing of high dimensional sparse functions \cite{r14}etc. It also has been employed in biomedical fields, i.e. developing the computational molecular modelling tools \cite{r15} and analyzing the level set estimation in medical image and pattern recognition \cite{r16} etc. Due to its widely applications, the refined results and improvements of Hoeffding's inequality and Hoeffding-Azuma inequality in martingales usually resulted in more new insights on the developments of related fields.  Recently, Hertz \cite{r1} presented an improvement result on the original Hoeffding's inequality by utilizing the asymmetric feature of finite distribution interval of random variables. It can reduce the related exponential coefficient from its arithmetic means to the geometric means of $|a|$ and $b$, where $[a,b]$ ($a<0, b>0)$ is the distributed interval of random variable $X$. This improvement motivate us to improve the Hoeffding's inequality.  For simplicity, let us first review the result of Hoeffding's inequality \cite{r2} and its improvement obtained by Hertz \cite{r1}.

\subsection{Hoeffding's Inequality and An Improvement}
Assume that $X$ is a zero mean real valued random variable and $X\in [a, b]$ with $a<0,b>0$. Hoeffding's lemma state that for all $s\in \textbf{R}, s>0$,
\begin{equation}\label{hoeffding inequality}
  E[e^{sX}]\leq \exp\Big\{\frac{s^2 (b-a)^2}{8}\Big\}
\end{equation}

Recently, D. Hertz presented an improved result with the following form
\begin{equation} \label{hertz ine}
 E[e^{sX}]\leq  \exp\Big\{\frac{s^2 \Phi^2(a,b)}{2}\Big \},
 \end{equation}
where
\begin{equation}
  \Phi(a,b)=\begin{cases} \frac{|a|+b}{2}\quad b >\mid a\mid, \\
  \sqrt{|a|b},\quad b\leq \mid a\mid.
  \end{cases}
\end{equation}

Since $\sqrt{|a|b} \leq \frac{|a|+b}{2}$, it gives a tighter upper bound for $-a>b$, compared with the original Hoeffding's inequality.

Motivated by this result, an interesting question is raised. Can we further improve the Hoeffding's inequality? If so, how to do it.

In this paper, we try to derive a new type of Hoeffding's inequalities, where higher order moments of random variable $X$ are taken into account, except $E(X)=0$. i.e. $E(X^k)=m_k (k=2,3,...)$.

\subsection{Main Theorem}
To give a clear picture of this paper, the new type of Hoeffding's inequalities are given as follows.

\begin{thm}\label{main theorem}
Assume that $X$ is a real valued random variable with $E(X)=0$, $X\in [a, b]$ with $a<0,b>0$. For all $s\in \textbf{R}, s>0$ and an integer $k$ ($k\geq 1$), then we have
\begin{equation}\label{hoeffding k}
  E[e^{sX}]\leq  \Upsilon_k(a,b) \exp\Big\{\frac{s^2}{2k}\Phi^2(a,b)\Big\}
\end{equation}
where
\begin{equation}
\Upsilon_k(a,b)=\Big[1+\frac{\max\{|a|,b\}}{|a|}\Big]^k-k\frac{\max\{|a|,b\}}{|a|}
\end{equation}
\begin{equation}
  \Phi(a,b)=\begin{cases} \frac{|a|+b}{2}\quad b >\mid a\mid, \\
  \sqrt{|a|b},\quad b\leq \mid a\mid.
  \end{cases}
\end{equation}
\end{thm}

\textbf{Remark 1.} When $k=1$, it is easy to check that $\Upsilon_1(a,b)=1$. This indicates that the new type Hoeffding's inequality will be reduced to the improved Hoeffding's inequality (\ref{hertz ine}),still better than the original Hoeffding's inequality. When $k=2$, $\Upsilon_1(a,b)=1+\{\frac{\max\{|a|,b\}}{|a|}\}^2$ and the exponential coefficient has been decreased by 2 times compared to the improved Hoeffding's inequality (\ref{hertz ine}). In fact, such a result can be refined, which is given by the following Corollary.

\begin{cor}\label{cor k equals 2}
Under the same assumption of Theorem \ref{main theorem}, for $k=2$, we have
\begin{equation}
  E[e^{sX}]\leq [1+\frac{m_2}{a^2}] \exp\Big\{\frac{s^2}{4}\Phi^2(a,b)\Big\}
\end{equation}
where $m_2=E(X^2)$.

If $E(X^2)$ is unknown, the inequalities can be relaxed as
\begin{equation}
  E[e^{sX}]\leq  [1+\frac{b}{|a|}] \exp\Big\{\frac{s^2}{4}\Phi^2(a,b)\Big\}
      \end{equation}
      and
    \begin{equation} \label{ine2}
  E[e^{sX}]\leq  2 \exp\Big\{\frac{s^2}{4}\Phi^2(a,b)\Big \} \quad \text{if } |a| \geq b
      \end{equation}
\end{cor}

Comparing the result in eqn.(\ref{ine2}) with that presented in Theorem \ref{main theorem}, it is easy to check that
 \begin{equation}
  [1+\frac{b}{|a|}]\leq 1+\{\frac{\max\{|a|,b\}}{|a|}\}^2
 \end{equation}
 holds.  This indicates that Corollary \ref{cor k equals 2} really improves the result presented in Theorem \ref{main theorem} for $k=2$. Comparing to that in eqn. (\ref{hertz ine}),  the exponential coefficient has be reduced by  2 times. That is to say, when parameter $s$ is relatively large, the new type of Hoeffding's inequalities will give much tighter results than original Hoeffding's inequality and its improvement obtained by Hertz.

The remaining part of this paper is organized as follows. In Section 2, we first present the  proof of Corollary (\ref{cor k equals 2}) and show the insight by taking higher order moments of real valued random variables into account and then present the proof of main theorem in this paper. In Section 3, we present the new type Hoeffding's inequalities applications in the one sided and two sided tail bound. We also discuss how to select the integer parameter $k$ to give a tighter bound In Section 4. Finally, in Section 5, we give the conclusion.

\section{The Proof of Main Theoretical Results}

Let us first introduce some Lemmas.

\subsection{Some Useful Lemmas}
\begin{lem}\label{lemma1}
Supposed $f(x)$ is a convex function of $x$, $f(x)>0$ with $x\in [a,b]$, then we have the following results.
(i)
\begin{equation*}
  f(x)\leq \frac{b-x}{b-a}f(a)+\frac{x-a}{b-a}f(b)
\end{equation*}
(ii) $f^2(x)$ is also a convex function of $x$ and
\begin{equation*}
  f^2(x)\leq [\frac{b-x}{b-a}f(a)+\frac{x-a}{b-a}f(b)]^2
\end{equation*}
and
\begin{equation*}
  f^2(x)\leq \frac{b-x}{b-a}f^2(a)+\frac{x-a}{b-a}f^2(b)
\end{equation*}
\end{lem}
The proof of Lemma \ref{lemma1} can be directly derived by using the definition of Convex function and $(f^2(x))'=2f(x)f'(x)$ and $(f^2(X))''=2(f'(x))^2+2f(x)f''(x)>0$ .

\begin{lem}\label{lemma2}
 Assume that $X$ is a real valued random variable with $E(X)=0$, $P(X\in [a, b])=1$ with $a<0,b>0$. We have
(i) \begin{equation}\label{x2}
  E(X^2)\leq |a|b
\end{equation}
(ii) \begin{equation}\label{x4}
  E(X^4)\leq |a|b(a^2+ab+b^2)\leq|a|b(a^2+b^2)
\end{equation}
\end{lem}

\begin{IEEEproof}
(i) Since $f(x)=x^2$ is a convex function of $x$ in $[a,b]$,  we have
\begin{equation}
  x^2\leq\frac{b-x}{b-a}a^2+\frac{x-a}{b-a}b^2
\end{equation}

\begin{equation}
\begin{split}
  E(X^2){} &\leq \frac{b}{b-a}a^2+\frac{-a}{b-a}b^2{}\\
  & =|a|b
\end{split}
\end{equation}
(ii)
Since $f(x)=x^2$ is a convex function and $f(x)\geq 0$, we know that $f^2(x)=x^4$ is also a convex function of $x$ according to Lemma \ref{lemma1}. Then we have
\begin{equation}\label{x4 estimate}
  x^4\leq\frac{b-x}{b-a}a^4+\frac{x-a}{b-a}b^4
\end{equation}
\begin{equation}
\begin{split}
  E(X^4){} &\leq \frac{b}{b-a}a^4+\frac{-a}{b-a}b^4{}\\
 & =|a|b (a^2+ab+b^2)\\
 & \leq |a| b(a^2+b^2)
\end{split}
\end{equation}
\end{IEEEproof}

\begin{lem}\label{lemma3}
For $0<\lambda<1$ and $u>0$, let
\begin{equation}
  \psi(u)=-\lambda u+\ln(1-\lambda+\lambda e^u)
\end{equation}
Then we have
\begin{equation}
  \psi (u) =0.5\tau (1-\tau)u^2
\end{equation}
where
\begin{equation}
  \tau =\frac{\lambda}{(1-\lambda)e^{-\xi}+\lambda}, \quad \xi\in[0,u]
\end{equation}
In addition, we have
\begin{equation}
\psi (u) \leq \begin{cases} \frac{u^2}{8}\quad \quad \quad \quad \quad \lambda \leq 0.5 \\
\lambda (1-\lambda) \frac{u^2}{2} \quad \lambda >0.5
\end{cases}
\end{equation}
\end{lem}
This lemma was derived in \cite{r1}. For completeness, we reorganize it as follows.
\begin{IEEEproof}

Since
\begin{equation}
  \psi(u)=-\lambda u+\ln(1-\lambda+\lambda e^u)
\end{equation}
For $u>0$, one can use Taylor's expansion and obtain
\begin{equation}
  \psi(u)=\psi(0)+\psi^{'}(0)u+0.5\psi^{''}(\xi)u^2
\end{equation}
it is easy to check that $\psi(0)=0$ and
\begin{equation}
  \psi^{'}(u)=-\lambda + \frac{\lambda e^u}{1-\lambda +\lambda e^u}
\end{equation}
\begin{equation}
  \psi^{''}(u)= \frac{\lambda e^u}{1-\lambda +\lambda e^u}(1-\frac{\lambda e^u}{1-\lambda +\lambda e^u})
\end{equation}
 That means
 \begin{equation*}
  \psi^{'}(0)=0
 \end{equation*}
  and
 \begin{equation}
  \psi^{''}(\xi)= 0.5 \tau (1-\tau)
\end{equation}
where
\begin{equation}
\tau =\frac{\lambda}{(1-\lambda)e^{-\xi}+\lambda}, \quad \xi\in[0,u] .
\end{equation}
That is,
\begin{equation}
\psi(u)=0.5 \tau (1-\tau)u^2
\end{equation}
Now let us divide it into two cases to discuss.
(a) If $\lambda >0.5$, then
\begin{equation}
   \tau=\frac{\lambda}{(1-\lambda)e^{-\xi}+\lambda}\geq\lambda>0.5
 \end{equation}
 That means,
  $\tau(1-\tau)$ reaches its maximum at $\tau=\lambda$,
  In other word,
  $\tau(1-\tau)\leq \lambda (1-\lambda)$
 (b) If $\lambda \leq 0.5$, then we have $\tau(1-\tau)\leq \frac{1}{4}$.

  By combining cases (a) and (b), we get
  \begin{equation}
\psi (u) \leq \begin{cases} \frac{u^2}{8}\quad \quad \quad \quad \quad \lambda \leq 0.5 \\
\lambda (1-\lambda) \frac{u^2}{2} \quad \lambda >0.5
\end{cases}
\end{equation}
The proof is completed.

\end{IEEEproof}

\subsection{Observation from Corollary \ref{cor k equals 2}}

Now we first review the Corollary \ref{cor k equals 2}. It claimed that under the same assumption of Theorem \ref{main theorem}, for $k=2$, we have
\begin{equation}
  E[e^{sX}]\leq [1+\frac{m_2}{a^2}] \exp\Big\{\frac{s^2}{4}\Phi^2(a,b)\Big\}
\end{equation}
where $m_2=E(X^2)$.

Before we present the proof of Corollary \ref{cor k equals 2}, let us analyze why such a new type of Hoeffding's inequality can decrease its exponential factor by 2 times in philosophy.
Since
\begin{equation*}
f(x)=\exp(\alpha x)
\end{equation*}
 is a convex function for any $\alpha>0$. \\

Let $\alpha =2\tilde{s}$, then
\begin{equation}
\begin{split}
  E(\exp(2\tilde{s}X))\leq {} &\frac{b^2+m_2}{(b-a)^2}\exp(2\tilde{s}a)+\frac{m_2+a^2}{(b-a)^2}\exp(2\tilde{s}b){}\\
  &+\frac{-2ab-2m_2}{(b-a)^2}\exp(\tilde{s}a)\exp(\tilde{s}b)\\
  & =\frac{b^2+m_2}{(b-a)^2}\exp(2\tilde{s}a)+\frac{m_2+a^2}{(b-a)^2}\exp(2\tilde{s}b)\\
 & +\frac{-2ab-2m_2}{(b-a)^2}\exp(2\tilde{s}\frac{a+b}{2})
  \end{split}
\end{equation}

The equation above can be rewritten as
\begin{equation}\label{alpha 2}
\begin{split}
  E(\exp(\alpha X))\leq &{} \frac{b^2+m_2}{(b-a)^2}\exp(\alpha a)+\frac{m_2+a^2}{(b-a)^2}\exp(\alpha b){}\\
  &+\frac{-2ab-2m_2}{(b-a)^2}\exp(\alpha \frac{a+b}{2})
  \end{split}
\end{equation}

Using Lemma \ref{lemma2} above, it is easy to see that all of the weighting coefficients of $\exp(\alpha a),\exp(\alpha b)$ and $\exp(\alpha \frac{a+b}{2})$
are non-negative and
\begin{equation}
  \frac{b^2+m_2}{(b-a)^2}+\frac{m_2+a^2}{(b-a)^2}+\frac{-2ab-2m_2}{(b-a)^2}=1
\end{equation}
Now, by using $s=\alpha$ in the inequality (\ref{alpha 2}), we have
\begin{equation} \label{s 2}
\begin{split}
  E(\exp(s X))\leq {}&\frac{b^2+m_2}{(b-a)^2}\exp(sa)+\frac{m_2+a^2}{(b-a)^2}\exp(sb){}\\
  &+\frac{-2ab-2m_2}{(b-a)^2}\exp(s\frac{a+b}{2})
  \end{split}
\end{equation}
It is easy to see that the right hand side of equation is equal to the linear weighting sum of $\exp(sa),\exp(sb)$ and $\exp(s\frac{a+b}{2})$. That is to say, one can use the information provided by three points to estimate the upper bound of $E(\exp(sX))$. It exactly provides more information than that only using two point linear weighting sum of $\exp(sa)$ and $\exp(sb)$ to estimate the upper bound of $E(\exp(sX))$.
Similarly, if one can use the information of function $\exp(sx)$ at multiple points, the upper bound of estimation $E(\exp(sX))$ may be improved further, this is why we consider the high order moments of random variables to discuss Hoeffding's inequality improvements.

\subsection{Proof of Corollary \ref{cor k equals 2}}

Now let us present the proof of Corollary \ref{cor k equals 2}.

\begin{IEEEproof}

Following the inequality (\ref{s 2}), we have that
\begin{equation}\label{s 2a}
\begin{split}
  E(\exp(s X))\leq {}&\frac{b^2+m_2}{(b-a)^2}\exp(sa)+\frac{m_2+a^2}{(b-a)^2}\exp(sb){}\\
  &+\frac{-2ab-2m_2}{(b-a)^2}\exp(s\frac{a+b}{2})\\
 & =\frac{b^2}{(b-a)^2}\exp(sa)+\frac{a^2}{(b-a)^2}\exp(sb)\\
  &+\frac{-2ab}{(b-a)^2}\exp(s\frac{a+b}{2})\\
 & +\frac{m_2}{(b-a)^2} \{ \exp(s\frac{b}{2})-\exp(s\frac{a}{2})\}
  \end{split}
\end{equation}

Let $u=\frac{s(b-a)}{2}$,$\lambda=\frac{-a}{b-a}$, and $\beta^2=\frac{m_2}{(b-a)^2}$, then we have
$s=\frac{2u}{b-a}$, $\frac{b}{b-a}=1-\lambda$. The inequality \ref{s 2a} can be rewritten as
\begin{equation}
  \begin{split}
  E(\exp(sX))\leq{}& [(1-\lambda)e^{-\lambda u}+\lambda e^{(1-\lambda)u}]^2 +\beta^2 (e^{(1-\lambda)u}-e^{-\lambda u})^2{} \\
  & \leq [(1-\lambda)e^{-\lambda u}+\lambda e^{(1-\lambda)u}]^2 (1+\frac{\beta^2}{\lambda^2})\\
  &=\exp(2\psi(u))(1+\frac{\beta^2}{\lambda^2})
  \end{split}
\end{equation}

By using Lemma \ref{lemma3}
 we have
 \begin{equation}
    E(\exp(sX))\leq \begin{cases} \exp\big( \frac{u^2}{4}\big)(1+\frac{\beta^2}{\lambda^2}) \quad \quad \quad \quad  |a|<b \\
      \exp(\lambda(1-\lambda)u^2)(1+\frac{\beta^2}{\lambda^2}) \quad |a|\geq b
      \end{cases}
     \end{equation}

 Now we shall discuss the exponential coefficient and the multiple factor $\big(1+\frac{\beta^2}{\lambda^2}\big)$ in tow different cases.

 (a)\quad If $|a|\geq b$, then by using $u=\frac{s(b-a)}{2}$,$\lambda=\frac{-a}{b-a}$, and $\beta^2=\frac{m_2}{(b-a)^2}$, we have
 \begin{equation}
  \lambda(1-\lambda)u^2\leq \frac{-ab}{(b-a)^2}\frac{s^2(b-a)^2}{4}=\frac{s^2|a|b}{4}
 \end{equation}
 as well as
 \begin{equation}
  1+\frac{\beta^2}{\lambda^2}=1+\frac{m_2}{a^2}\leq 1+\frac{b}{|a|}\leq 2
 \end{equation}

(b)\quad If $|a|<b$, then we have
\begin{equation}
  \frac{u^2}{4}=\frac{1}{4}\frac{s^2(b-a)^2}{4}=\frac{s^2(b-a)^2}{16}
 \end{equation}
 and
 \begin{equation}
  1+\frac{\beta^2}{\lambda^2}=1+\frac{m_2}{a^2}\leq 1+\frac{b}{|a|}
 \end{equation}

Combining the two difference cases and using
 \begin{equation*}
  \Phi(a,b)=\begin{cases} \frac{|a|+b}{2}\quad b >\mid a\mid, \\
  \sqrt{|a|b},\quad b\leq \mid a\mid.
  \end{cases}
\end{equation*}
 we obtain the result in inequality. The proof is completed.
\end{IEEEproof}

\subsection{Proof of Theorem \ref{main theorem}}

\begin{IEEEproof}
If $k=1$, it is the improved Hoeffding's inequality \ref{hertz ine}.

Now we mainly focus on the case of $k\geq 2$.

Since $f(x)=e^{\alpha x}$ is a convex function of $x$ for all $\alpha >0$ and $f(X)>0$, we have
\begin{equation}
  e^{\alpha x}\leq \frac{b-x}{b-a}e^{\alpha a}+\frac{x-a}{b-a}e^{\alpha b}
\end{equation}
and

for an positive integer $k$ ($k\geq2$), we have
\begin{equation}
\begin{split}
   e^{k\alpha x}{}&\leq \Big[\frac{b-x}{b-a}e^{\alpha a}+\frac{x-a}{b-a}e^{\alpha b}\Big]^k{}\\
   &=\Big\{\Big[\frac{b}{b-a}e^{\alpha a}+\frac{-a}{b-a}e^{\alpha b}\Big]+x\Big[\frac{e^{\alpha b}-e^{\alpha a}}{b-a}\Big]\Big\}^k
      \end{split}
\end{equation}
and
\begin{equation}
   E\big(e^{k\alpha X}\big)\leq \Big\{\Big[\frac{b}{b-a}e^{\alpha a}+\frac{-a}{b-a}e^{\alpha b}\Big]+X\Big[\frac{e^{\alpha b}-e^{\alpha a}}{b-a}\Big]\Big\}^k\\
   \end{equation}
   By using $ s=k\alpha$ and $\lambda=\frac{-a}{b-a}$, $u=\frac{s}{k}(b-a)$, then we have
\begin{equation}
\begin{split}
   E\big(e^{s X}\big){}&\leq E\Big\{[(1-\lambda)e^{-\lambda u}+\lambda e^{(1-\lambda)u}]{}\\
   &+\frac{X}{|a|}[\lambda e^{(1-\lambda)u}-\lambda e^{-\lambda u}]\Big\}^k
   \end{split}
\end{equation}
Let
\begin{equation}
e^{\psi(u)}=(1-\lambda)e^{-\lambda u}+\lambda e^{(1-\lambda)u}
\end{equation}

 and
 \begin{equation}
 \varphi(u)=\lambda e^{(1-\lambda)u}-\lambda e^{-\lambda u}
 \end{equation}

then
   \begin{equation}
   \begin{split}
    E\big(e^{s X}\big)\leq{}& E\Big[e^{\psi(u)}+\frac{X}{|a|}\varphi(u)\Big]^k\\
    &= e^{k\psi(u)}+ke^{(k-1)\psi(u)}E\Big(\frac{X}{|a|}\Big)\varphi(u)\\
    &+\sum_{i=2}^k \mathcal{C}_k^i e^{(k-i)\psi(u)}E\Big(\frac{X}{|a|}\Big)^i\varphi^i(u)
    \end{split}
    \end{equation}
    \begin{equation}
    \begin{split}
    =&e^{k\psi(u)}+\sum_{i=2}^k \mathcal{C}_k^i e^{(k-i)\psi(u)}E\Big(\frac{X}{|a|}\Big)^i\varphi^i(u)\\
    &\leq e^{k\psi(u)}+\sum_{i=2}^k \mathcal{C}_k^i e^{(k-i)\psi(u)}E\Big(\frac{|X|}{|a|}\Big)^i\varphi^i(u)\\
        &  \leq [(1-\lambda)e^{-\lambda u}+\lambda e^{(1-\lambda)u}]^k\\
       & \times \Big\{\Big[1+\frac{\max\{-a,b\}}{|a|}\Big]^k-k\frac{\max\{-a,b\}}{|a|}\Big\}
       \end{split}
   \end{equation}
   where $\mathcal{C}_k^i=\frac{k!}{i!(k-i)!}$.

By using $(b-a)\lambda=-a$, and $\psi(u)=0.5 \tau (1-\tau)u^2$, $u=\frac{s}{k}(b-a)$, we have
\begin{equation}
\begin{split}
  E(e^{s X})&\leq e^{\frac{k}{2}\tau(1-\tau)u^2} \Big\{[1+\frac{\max\{-a,b\}}{|a|}]^k-k\frac{\max\{-a,b\}}{|a|}\Big\} \\
  &\leq \Big \{[1+\frac{\max\{|a|,b\}}{|a|}]^k-k\frac{\max\{-a,b\}}{|a|}\Big \} \exp\Big(\frac{s^2}{2k}\Phi^2\Big)\\
  &= \Upsilon_k (a,b)\exp\Big(\frac{s^2}{2k}\Phi^2\Big)
  \end{split}
\end{equation}
where
\begin{equation}
\Upsilon_k(a,b)=\Big[1+\frac{\max\{|a|,b\}}{|a|}\Big]^k-k\frac{\max\{|a|,b\}}{|a|},
\end{equation}
\begin{equation}
  \Phi = \begin{cases} \frac{(b-a)}{2} \quad  -a<b \\
      \sqrt{|a|b} \quad -a \geq b
  \end{cases}
\end{equation}
The proof is completed.
\end{IEEEproof}

\textbf{Remark 2.}
The proof of theorem \ref{main theorem} create a new routine on how to use multipoint values of $\exp{sx}$ to get tighter approximation of $E(\exp{sX})$ for any random distribution in a finite interval with $P(X\in [a,b])=1$.  Comparing with the original Hoeffding's inequality and its improvement obtained by Hertz,  the advantages is  that it can exactly reduce the exponential coefficients by $k$ times when all the moments of less than $k$ order statistics are taken into account, but the cost is that it will almost enlarge the multiply factor with $C_1^k$ times, as shown by $\Upsilon_k(a,b)$, where $C_1$ is a constant with $C_1>1$ . That means there exists a trade off between the exponential coefficient reduction and the multiply factor increment. It needs to be consider in specific applications.

In some scenarios, one may interest in the case of $k=4$. The following Corollary shows one refinement of Theorem \ref{main theorem}.

\begin{cor}\label{Cor k equal 4}
Assume that $X$ is a real valued random variable, $P(X\in [a, b])=1$ with $a<0, b>0$ and  $E(X)=0$, $E(X^2)=m_2$, $E(X^3)=0$ and $E(X^4)=m_4$.  For all $s\in \textbf{R}, s>0$, we have
\begin{equation}\label{hoeffding inequality3}
  E[e^{sX}]\leq \Big[1+\frac{6m_2}{a^2}+\frac{m_4}{a^4}\Big] \exp\Big(\frac{s^2}{8}\Phi^2(a,b)\Big)
\end{equation}
where \begin{equation}
  \Phi = \begin{cases} \frac{(b-a)}{2} \quad  -a<b \\
      \sqrt{|a|b} \quad -a \geq b
  \end{cases}
\end{equation}
\end{cor}

\begin{IEEEproof}
The proof can follow that way on Theorem\ref{main theorem}.

Since $f(x)=e^{\alpha x}$ is a convex function of $x$ for all $\alpha >0$ and $f(X)>0$, we have
\begin{equation}
  e^{\alpha x}\leq \frac{b-x}{b-a}e^{\alpha a}+\frac{x-a}{b-a}e^{\alpha b}
\end{equation}
and
\begin{equation}
\begin{split}
   e^{4\alpha x} {}&\leq \Big[\frac{b-x}{b-a}e^{\alpha a}+\frac{x-a}{b-a}e^{\alpha b}\Big]^4{}\\
   &=\Big\{\Big[\frac{b}{b-a}e^{\alpha a}+\frac{-a}{b-a}e^{\alpha b}\Big]+x\Big[\frac{e^{\alpha b}-e^{\alpha a}}{b-a}\Big]\Big\}^4
      \end{split}
\end{equation}
Let $s=4\alpha$, and using  $E(X)=0$,$E(X^2)=m_2$ $E(X^3)=0$ and $E(X^4)=m_4$, we have
\begin{equation}
\begin{split}
  E(e^{sX})\leq {}&\Big(\frac{b}{b-a}e^{\frac{s}{4} a}+\frac{-a}{b-a}e^{\frac{s}{4} b}\Big)^4{}\\
  &+6m_2\Big(\frac{b}{b-a}e^{\frac{s}{4} a}+\frac{-a}{b-a}e^{\frac{s}{4} b}\Big)^2\Big(\frac{e^{\frac{s}{4} b}-e^{\frac{s}{4} a}}{b-a}\Big)^2\\
  &+m_4\Big(\frac{e^{\frac{s}{4} b}-e^{\frac{s}{4} a}}{b-a}\Big)^4
  \end{split}
\end{equation}
\\
Let $\lambda=\frac{-a}{b}$, $u=\frac{s}{4}(b-a)$, then we have
$\frac{b}{b-a}=1-\lambda$, $ \frac{s}{4}a=-\lambda u$, $\frac{s}{4}b=(1-\lambda) u$.
Then the inequality above can be rewritten as
\begin{align}
  E(e^{sX})\leq & [(1- \lambda) e^{-\lambda u}+\lambda e^{(1-\lambda)u}]^4 \\
  &+\frac{6 m_2}{(b-a)^2\lambda^2}[(1-\lambda)e^{-\lambda u}+\lambda e^{(1-\lambda)u}]^2\\
   & \times [\lambda e^{(1-\lambda)u}-\lambda e^{-\lambda u}]^2
   \end{align}

 \begin{equation}
 \begin{split}
  & +\frac{m_4}{(b-a)^4\lambda^4}[\lambda e^{(1-\lambda)u}-\lambda e^{-\lambda u}]^4\\
  \leq &[(1-\lambda) e^{-\lambda u}+\lambda e^{(1-\lambda)u}]^4\\
  &\times [1+\frac{6m_2}{(b-a)^2\lambda^2}+\frac{m_4}{(b-a)^4\lambda^4}]\\
  =& e^{4\psi(u)}[1+\frac{6m_2}{(b-a)^2\lambda^2}+\frac{m_4}{(b-a)^4\lambda^4}]
    \end{split}
\end{equation}
by using $(b-a)\lambda=-a$, and Lemma \ref{lemma3}, we have
\begin{equation}
   E(e^{sX})\leq [1+\frac{6m_2}{a^2}+\frac{m_4}{a^4}]e^{\frac{s^2\Phi^2}{8}}
\end{equation}
where
\begin{equation}
  \Phi = \begin{cases} \frac{(b-a)}{2} \quad  -a<b \\
      \sqrt{|a|b}\quad -a \geq b
  \end{cases}
\end{equation}

The proof is completed.
\end{IEEEproof}
If the $E(X^2)$ and $E(X^4)$ are not exactly known and $|a|=b$, we have the following result.

\begin{cor}\label{cor k equal 4 and a=b}
Assume that $X$ is a real valued random variable.$P(X\in [-a, a])=1$ with $a>0$ and  $E(X)=0$ and $E(X^3)=0$ .  For all $s\in \textbf{R}, s>0$, we have
\begin{equation}
  E[e^{sX}]\leq 8 \exp(\frac{a^2s^2}{8})
  \end{equation}
\end{cor}

\begin{IEEEproof}
by using
$m_2\leq a^2$ and $m_4 \leq a^4$ and the inequality in Corollary \ref{Cor k equal 4} , we can get the result directly.
\end{IEEEproof}

\section{Applications in Tail Bound Estimation}
Let us consider the scenario, where $X_1,X_2,\dots,X_n$ be independent random variables such that $X_i\in [a_i, b_i], a_i<0, b_i>0$ and $EX_i=0$ for $i=1,2,\dots,n$. Define $S_n=\sum_{i=1}^n$.
It is easy to check that $ES_n=0$.  For all $s>0$, we have
\begin{equation}
 \begin{split}P(S_n\geq t){}& = P\big(e^{sS_n}\geq e^{st}\big) \quad \text{Chernoff} {}\\
&\leq e^{-st}Ee^{sS_n} \quad \quad \quad \text{Markov}\\
&= e^{-st} \prod _{i=1}^n Ee^{sX_i}
\end{split}
\end{equation}

Using the results of Theorem \ref{main theorem} and its Corollaries, one can obtain that
\begin{equation}
Ee^{sX_i}\leq A_{k_i} \exp\Big(\frac{s^2}{2k_i}\Phi_i\Big)
\end{equation}
where $A_{k_i}$ and $k_i$ are based on which one inequality of $X_i$ being selected in Section II and III with
\begin{equation}
A_{k_i}= \begin{cases} 1 \quad \quad \quad \quad \quad \quad \quad \quad \quad \quad \quad \quad \quad \quad k_i=1\\
1+\frac{\max\{|a|,b\}}{|a|} \quad \quad \quad \quad \quad \quad \quad \quad  \quad k_i=2\\
\Big[1+\frac{\max\{|a|,b\}}{|a|}\Big]^{k_i}-k_i\frac{\max\{|a|,b\}}{|a|} \quad k_i\geq 3
\end{cases}
\end{equation}
and
$\Phi_i=\Phi(a_i,b_i)$.

In this case,  we get
\begin{equation}
P(S_n\geq t)\leq \Big(\prod_{i=1}^n A_{k_i}\Big) \exp\Big\{-st+s^2\Big(\sum_{i=1}^n \frac{\Phi_i^2}{2k_i}\Big)\Big\}
\end{equation}

Now selecting
\begin{equation}
s=\frac{t}{2\Big(\sum_{i=1}^n \frac{\Phi_i^2}{2k_i}\Big)}
\end{equation}
to minimize the exponent in the last inequality, we obtain
\begin{equation} \label{one side tail bound}
P(S_n\geq t)\leq \Big(\prod_{i=1}^n A_{k_i}\Big)\exp\Big\{-t^2\Big(2\sum_{i=1}^n \frac{\Phi_i^2}{k_i}\Big)^{-1}\Big\}
\end{equation}
In particular, if all the $k_i$, $(i=1,2,\dots,n)$ are selected as 1, then $A_{k_i}=1$, it reduces to the Improved Hoeffding's one side tail bound.

If all the $k_i$, $(i=1,2,\dots,n)$ are selected as 2 and $|a_i|=b_i$, then $A_{k_i}=2$, and the inequality can be rewritten as
\begin{equation}
P(S_n\geq t)\leq 2^n \exp\Big\{-\frac{t^2}{\sum_{i=1}^n a_i^2}\Big\}
\end{equation}

Furthermore,
\begin{equation}
P\Big(\frac{S_n}{n}\geq l \Big)\leq \Big(\prod_{i=1}^n A_{k_i}\Big) \exp\Big\{\frac{-nl^2}{2\tilde{\Phi}_i^2}\Big\}
\end{equation}
where $l$ is a positive number and
$\tilde{\Phi}_i^2=\frac{1}{n}\Big(\sum_{i=1}^n \frac{\Phi_i^2}{2k_i}\Big)$.

The two sided tail bound can be given by
\begin{equation}
\begin{split}
P(|S_n|\geq t){}&\leq \Big(\prod_{i=1}^n A_{k_i}\Big)\exp\Big\{-t^2\Big(2\sum_{i=1}^n \frac{\Phi_i^2}{k_i}\Big)^{-1}\Big\}{}\\
&+\Big(\prod_{j=1}^n B_{k_j}\Big)\exp\Big\{-t^2\Big(2\sum_{j=1}^n \frac{\Phi_j^2}{k_j}\Big)^{-1}\Big\}
\end{split}
\end{equation}
where $\{B_{k_j},j=1,2,\cdots, n\}$ is a sort of $\{A_{k_i},i=1,2,\cdots, n\}$ complement. That is to say, The calculation of $B_{k_j}$ is  just changing the positions of $a_j$ and $b_j$ in such a way $-a_j\rightarrow b_i$ and $-b_j\rightarrow a_i$ in the calculation of $A_{k_i}$ if the integer index $k_j$ of $B_{k_j}$ is equal to the integer index $k_i$ of  $A_{k_i}$.  In other word, for the same $X_i$,  it may select two different integer parameter values of $k_i$  to estimate both sided tail bounds for the positive and the negative directions.


\section{Selection of Integer Parameter $k$}

On the selection of integer parameter $k_i$, we shall discuss it firstly from one sided tail bound.
For simplicity, let us consider $n=1$. The first question is when selecting a larger $k$ will get a tighter bound.
The question can be solved by
\begin{equation}\label{comarison with k and k+1}
A_{k+1}\exp\Big\{-t^2\Big(2\frac{\Phi^2}{k+1}\Big)^{-1}\Big\}< A_{k}\exp\Big\{-t^2\Big(2\frac{\Phi^2}{k}\Big)^{-1}\Big\}
\end{equation}
Using logarithm on both sides of inequality (\ref{comarison with k and k+1}) and after some manipulations,  we get
\begin{equation}
\frac{t^2}{2\Phi^2}> \ln A_{k+1}-\ln A_k
\end{equation}
 That is
 \begin{equation}
 t > \Phi \sqrt{2 \ln \frac{A_{k+1}}{A_k}}
 \end{equation}

 To clear illustrate the effect of $k$ selection, we give three examples.

\emph{Example 1.}  $a= -1$, $b=1$.  The selection  rule of $k$ ($k=1,2,3$) is given by
\begin{equation}
k=
\begin{cases}
{}&1 \quad 0<t< \sqrt{2\ln2}\approx 1.177{} \\
 & 2 \quad  \sqrt{2\ln2} <t< \sqrt{2\ln (2.5)}\approx 1.3537\\
 & 3 \quad  t > \sqrt{2\ln (2.5)}
\end{cases}
\end{equation}

 \emph{Example 2.}  $a= -1$, $b=5$.  The selection  rule of $k$ ($k=1,2,3$) is given by
\begin{equation}
k=
\begin{cases}
{}&1 \quad 0<t< 3\sqrt{2\ln6}\approx 5.679{} \\
 & 2 \quad  3\sqrt{2\ln6} <t< 3\sqrt{2\ln (191/6)}\approx7.892 \\
 & 3 \quad  t > 3\sqrt{2\ln (191/6)}
\end{cases}
\end{equation}

 \emph{Example 3.} $a= -5$, $b=1$.  The selection  rule of $k$ ($k=1,2,3$) is given by
\begin{equation}
k=
\begin{cases}
{}&1 \quad 0<t< \frac{1}{2}\sqrt{10\ln(6/5)}\approx 0.6751 {} \\
 & 2 \quad  \frac{1}{2}\sqrt{10\ln(6/5)} <t< \sqrt{10\ln (25/6)}\approx 3.778 \\
 & 3 \quad  t > \sqrt{10\ln (25/6)}
\end{cases}
\end{equation}

\textbf{Remark 3.}
All the three examples shows that when $t$ is relatively small, i.e. close to zero, selecting parameter $k=1$ is the best one. The results in \emph{Example 3} shows that when $t=0.8$, selecting $k=2$ will give a tighter bias bound.  The results in \emph{Example 2} and \emph{Example 3} also indicates when random variable $X$ with $P(X\in [-1,5])=1$, where $a=-1, b=5$, one need to estimate $P(|X|> 0.8)$, the right hand sided bound should select $k=1$ as its estimation while the left hand sided bound should select $k=2$ as its estimation. This result shows that one may not consistently select the same parameter $k$  to deal with both sided bias bounds when  $|a|\neq b$.

Now let consider the general case.

The goal of parameters $k_i$ selection is basically to minimize the right hand of inequality (\ref{one side tail bound}). Thus, one can set up an optimization problem as follows.

\textbf{Problem 1}:
For a given $t>0$,
\begin{equation}
\min\limits_{k_i}\Big(\prod_{i=1}^n A_{k_i}\Big)\exp\Big\{-t^2\Big(2\sum_{i=1}^n \frac{\Phi_i^2}{k_i}\Big)^{-1}\Big\}
\end{equation}
where $A_{k_i}$ are calculated by using the theoretical results in Theorem \ref{main theorem} and its Corollaries for a given $k_i$.

It is equivalent to

\begin{equation}
\min\limits_{k_i}\Big(\sum_{i=1}^n \ln(A_{k_i})\Big)-t^2\Big(2\sum_{i=1}^n \frac{\Phi_i^2}{k_i}\Big)^{-1}
\end{equation}
and
\begin{equation}
\max\limits_{k_i}\frac{1}{\Big(2\sum_{i=1}^n \frac{\Phi_i^2}{k_i}\Big)}-\frac{\Big(\sum_{i=1}^n \ln(A_{k_i})\Big)}{t^2}
\end{equation}

In fact, such an optimization problem can be solved by using computer search.
Here, in order to provide a tractable mode, we relax $A_{k_i}$ with the form $\Big[1+\frac{\max\{|a|,b\}}{|a|}\Big]^k$ given in Theorem \ref{main theorem}.
In this case, the optimization problem can be transformed into the following problem.

\textbf{Problem 2: }
\begin{equation}
\max\limits_{k_i}\frac{1}{\Big(2\sum_{i=1}^n \frac{\Phi_i^2}{k_i}\Big)}-\frac{\Big(\sum_{i=1}^n k_i\ln\big(1+\frac{max\{|a_i|,b_i\}}{|a_i|}\big)\Big)}{t^2}
\end{equation}
Let us define
\begin{equation}
\begin{split}
g(k_1,k_2,\dots,k_n){}&=\frac{1}{\Big(2\sum_{i=1}^n \frac{\Phi_i^2}{k_i}\Big)}{}\\
&-\frac{\Big(\sum_{i=1}^n k_i\ln\big(1+\frac{max\{|a_i|,b_i\}}{|a_i|}\big)\Big)}{t^2}
\end{split}
\end{equation}

In order to get some insights, let us consider $k_j$ to be a real number rather than an integer. Then
the partial derivative of function $g(.)$ to $k_j$ is given by

\begin{equation}
\begin{split}
\frac{\partial g}{\partial k_j}= {}&\frac{2\Phi_i^2k_j^{-2}}{\Big(2\sum_{i=1}^n \frac{\Phi_i^2}{k_i}\Big)^2}{}\\
&-\frac{\ln\big(1+\frac{max\{|a_j|,b_j\}}{|a_j|}\big)}{t^2}
\end{split}
\end{equation}
Let $\frac{\partial g}{\partial k_j}=0$, after some manipulations, we obtain
\begin{equation}
\begin{split}
k_j={}&\frac{\Phi_j}{\sqrt{2\ln\big(1+\frac{max\{|a_j|,b_j\}}{|a_j|}\big) }}\frac{t}{\sum_{i=1}^n \frac{\Phi_i^2}{k_i}}{}\\
=& \frac{\Phi(a_j,b_j)}{\sqrt{2\ln\big(1+\frac{max\{|a_j|,b_j\}}{|a_j|}\big) }}\frac{t}{\sum_{i=1}^n \frac{\Phi_i^2}{k_i}}
\end{split}
\end{equation}
Since $ \frac{t}{\sum_{i=1}^n \frac{\Phi_i^2}{k_i}}$ is a common factor for all the $k_j$,$(j=1,2,\dots,n)$. This means
\begin{equation}
k_j \propto \frac{\Phi(a_j,b_j)}{\sqrt{2\ln\big(1+\frac{max\{|a_j|,b_j\}}{|a_j|}\big) }}
\end{equation}
That is to say, the near optimal value of $k_j$ is mainly determined by $a_j$ and $b_j$ except a common factor, the  parameters of distribution interval of  $X_j$.
This is an interesting result, which can provide more insight. In most of applications, all the $X_i$ $(i=1,2,\dots, n)$ are distributed with the same interval. In this case, one can select the same $k_i$ value for all of them, so that it can approximate the near optimal tighter tail bound.
Such a discussion can be extended to the scenarios of two sided tail bound.

\textbf{Remark 4.}
Consider the distribution interval is symmetric, where $|a_i|=b_i$. In this case, we have
\begin{equation}
\frac{\Phi(a_j,b_j)}{\sqrt{2\ln\big(1+\frac{max\{|a_j|,b_j\}}{|a_j|}\big) }}=\frac{|a_j|}{\sqrt{2\ln2}}
\end{equation}
This means
\begin{equation}
k_j \propto |a_j|
\end{equation}
This indicates the integer parameter $k_j$ selection is proportional to the distribution interval length. When $|a_j|$ is relatively small, i.e. $|a_j|$ is close to zero, the linear interpolation of two points with $x_1=a_j$ and $x_2=|a_j|$ is good enough to approximate the random curve of $e^{sX}$. That is to say, select $k_j=1$ is good enough.

When  $|a_j|$ is relatively large, the linear interpolation of two points with $x_1=a_j$ and $x_2=|a_j|$ may not be good enough to approximate the curve of $e^{sX}$. It needs more points in the curve of $e^{sx}$ to do the interpolation so that it could have a good approximation to the random curve of $e^{sX}$. That is to say, selecting a larger $k_j$ is necessary.  Such an observation is consistent with our "intuitive feeling" on the function approximation in philosophy.  We shall illustrate such phenomenon in detail with some examples below.

\emph{Example 4.} Let $a= -5$, $b=5$ and $m_2=5$.  The selection  rule of $k$ ($k=1,2,3$) is given by
\begin{equation}
k=
\begin{cases}
{}&1 \quad 0<t< 5\sqrt{2\ln2(6/5)}\approx 3.019 {} \\
 & 2 \quad  5\sqrt{2\ln2(6/5)} <t< 5\sqrt{2\ln (25/6)}\approx 8.447\\
 & 3 \quad  t > 5\sqrt{2\ln (25/6)}
\end{cases}
\end{equation}

\textbf{Remark 5.}
 The results in \emph{Example 1} shows selecting $k=1$ is always the best since the best working region for $t$  of $k \geq 2$ is out of the $X$ distributed interval, which can not occur in practice. \emph{Example 4} shows that when $m_2$ is given, it is possible to select $k\geq 2$ to get a tighter tail bound, i.e. $t=4$, the best selection of $k$ is $k=2$, which also show that when the distribution interval is relatively larger, it is possible to select the larger integer value of $k$ for the tail bound estimation.

\emph{Example 5.} Let us consider $n=4$, where $X_1 \in [-1,1]$, $X_2 \in [-5,5]$,  $X_3 \in [-1,5]$ and $X_4 \in [-5,1]$ with $E(X_1)=E(X_2)=E(X_3)=E(X_4)=0$. $ E(X_2^2)=5$ and $S_4=X_1+X_2+X_3+X_4$. It is easy to check that $S_4 \in [-12,12]$.
Fig. 1 shows different curves of one sided tail bounds , which are group one: $k_1=k_2=k_3=k_4 =1$. Group two: $k_1=k_3=k_4 =1, k_2=2$ and Group three: $k_1=k_3=1, k_2=k_4=2$, where the y-label is the logarithm of the one sided tail bound, $\Big(\sum_{i=1}^n \ln(A_{k_i})\Big)-t^2\Big(2\sum_{i=1}^n \frac{\Phi_i^2}{k_i}\Big)^{-1}$, the x-label is $t$.  It is observed that among the three groups of parameter $k$ selection, when $0<t<5.6647$, the curve of Group one provides the tightest bound. When $5.6647<t<10.0138$,  the curve of Group two provides the tightest bound and when $10.0138<t<12$,the curve of Group three provides the tightest bound.

\begin{figure}[http]
\centering
\includegraphics[width=220pt,height=220pt]{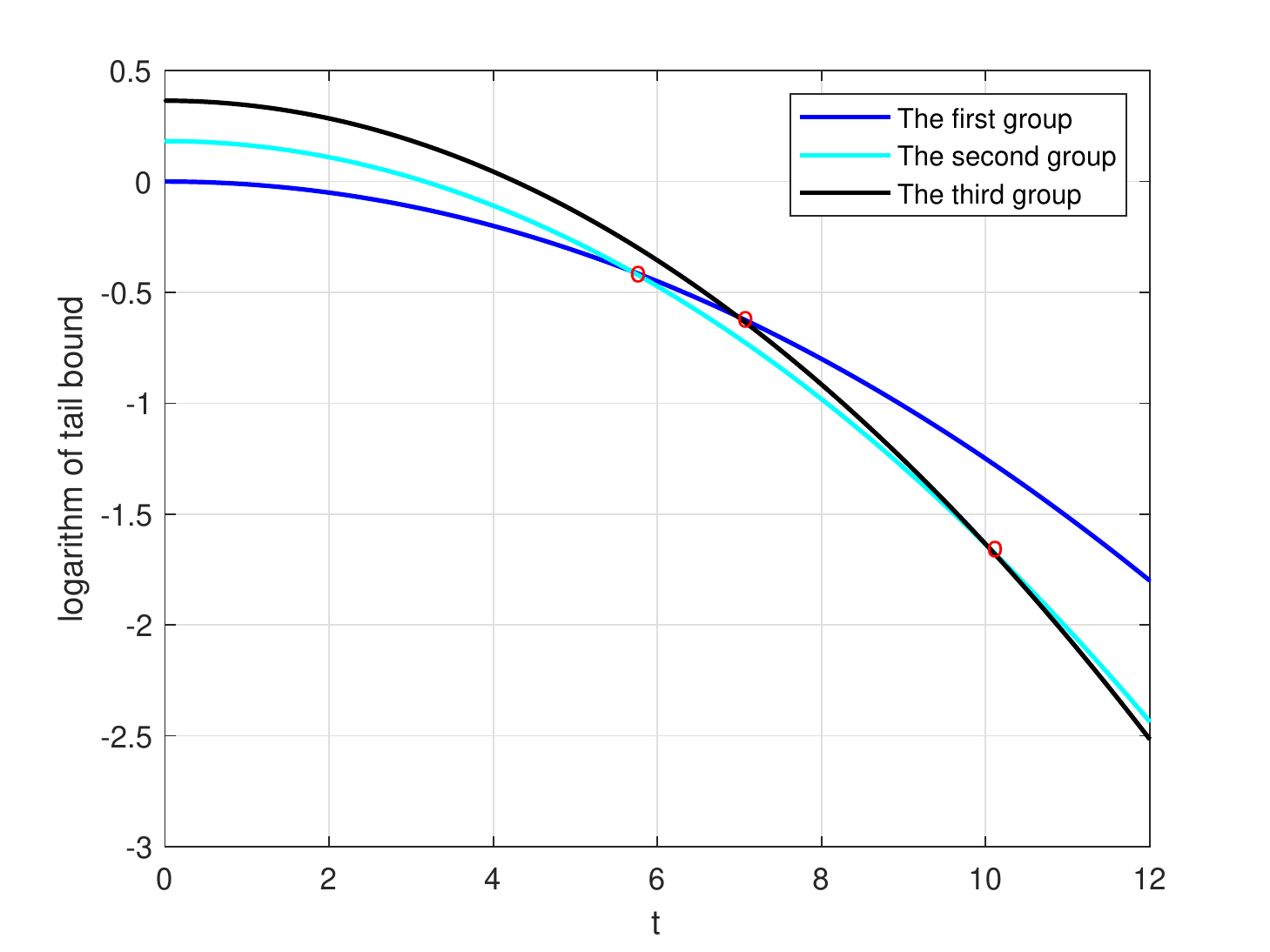}
 \caption{The logarithm of the one sided tail bound of three different group selection of parameters $k$ with $n=4$ in Example 5}
\label{fig:example 5}
\end{figure}


The results in \emph{Example 5} exactly demonstrated that the new type Hoeffding's inequalities are useful in the tail bound estimation.

\textbf{Remark 6.}
In real applications, one would not like to pay more attention on the selection of parameter $k_i$ in order to simply the system analysis. It recommends to select $k_i$ to be $ 1$ or $2$.

\section{Conclusion}

In this paper, we presented new type of Hoeffding's inequalities by using higher order moments of random variables. Some applications in one and two sided tail bound improvements can also be obtained by using the exponential function positiveness and Chernoff inequality. Perhaps, future research may focus on trying to improve the related inequalities that use Hoeffding's Lemma.


\end{document}